\documentclass{article}

\usepackage{graphicx,etoolbox}

\usepackage{hyperref}

\usepackage{hyperref}

\usepackage{amssymb}

\newcommand{\newc}{\newcommand}

\newc{\eqnoset}{\setcounter{equation}{0}}

\newcommand{\mref}[1]{(\ref{#1})}
\newcommand{\reflemm}[1]{Lemma~\ref{#1}}
\newcommand{\refrem}[1]{Remark~\ref{#1}}
\newcommand{\reftheo}[1]{Theorem~\ref{#1}}

\newcommand{\refsec}[1]{Section~\ref{#1}}

\newcommand{\beq}{\begin{equation}}
	\newcommand{\eeq}{\end{equation}}
\newcommand{\beqno}[1]{\begin{equation}\label{#1}}
	
	\newcommand{\barr}{\begin{array}}
		\newcommand{\earr}{\end{array}}
	
	\newc{\bearr}{\begin{eqnarray*}}
		\newc{\eearr}{\end{eqnarray*}}
	
	\newc{\bearrno}[1]{\begin{eqnarray}\label{#1}}
		\newc{\eearrno}{\end{eqnarray}}
	
	\newc{\non}{\nonumber}
	\newc{\nol}{\nonumber\nl}
	
	\newcommand{\bdes}{\begin{description}}
		\newcommand{\edes}{\end{description}}
	\newc{\benu}{\begin{enumerate}}
		\newc{\eenu}{\end{enumerate}}
	\newc{\btab}{\begin{tabular}}
		\newc{\etab}{\end{tabular}}

	\newtheorem{theorem}{Theorem}[section]
	\newtheorem{defi}[theorem]{Definition}
	\newtheorem{lemma}[theorem]{Lemma}
	\newtheorem{rem}[theorem]{Remark}
	\newtheorem{exam}[theorem]{Example}
	\newtheorem{propo}[theorem]{Proposition}
	\newtheorem{corol}[theorem]{Corollary}
	\newtheorem{conj}[theorem]{Conjecture}

	\newcommand{\btheo}[1]{\begin{theorem}\label{#1}}
		\newc{\brem}[1]{\begin{rem}\label{#1}\em}
			\newc{\bexam}[1]{\begin{exam}\label{#1}\em}
				\newc{\bdefi}[1]{\begin{defi}\label{#1}}
					\newcommand{\blemm}[1]{\begin{lemma}\label{#1}}
						\newcommand{\bprop}[1]{\begin{propo}\label{#1}}
							\newcommand{\bcoro}[1]{\begin{corol}\label{#1}}
								\newcommand{\btheoc}[1]{\begin{conj}\label{#1}}
									\newcommand{\etheo}{\end{theorem}}
								\newc{\etheoc}{\end{conj}}
							\newcommand{\elemm}{\end{lemma}}
						\newcommand{\eprop}{\end{propo}}
					\newcommand{\ecoro}{\end{corol}}
				\newc{\erem}{\end{rem}}
			\newc{\eexam}{\end{exam}}
		\newc{\edefi}{\end{defi}}
	
	\newc{\rmk}[1]{{\bf REMARK #1: }}
	\newc{\DN}[1]{{\bf DEFINITION #1: }}
	
	\newcommand{\bproof}{{\bf Proof:~~}}
	\newc{\eproof}{{\vrule height8pt width5pt depth0pt}\vspace{3mm}}
	
	\newc{\bfrac}[2]{\dspl{\frac{#1}{#2}}}

	\newc{\nid}{\noindent}
	
	\newcommand{\dspl}{\displaystyle}
	\newc{\grad}{\nabla}
	\newc{\Div}{\mbox{div}}
	\newc{\pdt}[1]{\dspl{\frac{\partial{#1}}{\partial t}}}
	\newc{\pdn}[1]{\dspl{\frac{\partial{#1}}{\partial \nu}}}
	\newc{\pdNi}[1]{\dspl{\frac{\partial{#1}}{\partial \mathcal{N}_i}}}
	\newc{\pD}[2]{\dspl{\frac{\partial{#1}}{\partial #2}}}
	\newc{\dt}{\dspl{\frac{d}{dt}}}
	\newc{\bdry}[1]{\mbox{$\partial #1$}}
	\newc{\sgn}{\mbox{sign}}
	
	\newc{\Hess}[1]{\frac{\partial^2 #1}{\pdh z_i \pdh z_j}}
	\newc{\hess}[1]{\partial^2 #1/\pdh z_i \pdh z_j}

	\newc{\ag}{\alpha}
	\newc{\bg}{\beta}
	\newc{\cg}{\gamma}\newc{\Cg}{\Gamma}
	\newc{\dg}{\delta}\newc{\Dg}{\Delta}
	\newc{\eg}{\varepsilon}
	\newc{\zg}{\zeta}
	\newc{\thg}{\theta}
	\newc{\llg}{\lambda}\newc{\LLg}{\Lambda}
	\newc{\kg}{\kappa}
	\newc{\rg}{\rho}
	\newc{\sg}{\sigma}\newc{\Sg}{\Sigma}
	\newc{\tg}{\tau}
	\newc{\fg}{\phi}\newc{\Fg}{\Phi}
	\newc{\vfg}{\varphi}
	\newc{\og}{\omega}\newc{\Og}{\Omega}
	\newc{\pdh}{\partial}
	
	\newc{\ccG}{{\cal G}}

	\newc{\ii}[1]{\int_{#1}}
	\newc{\iidx}[2]{{\dspl\int_{#1}~#2~dx}}
	\newc{\bii}[1]{{\dspl \ii{#1} }}
	\newc{\biii}[2]{{\dspl \iii{#1}{#2} }}
	\newc{\su}[2]{\sum_{#1}^{#2}}
	\newc{\bsu}[2]{{\dspl \su{#1}{#2} }}

	\newc{\biiom}[1]{{\dspl\int_{\bdrom}~ #1 ~d\sg}}
	\newc{\io}[1]{{\dspl\int_{\Og}~ #1 ~dx}}
	\newc{\bio}[1]{{\dspl\int_{\bdrom}~ #1 ~d\sg}}
	\newc{\bsir}{\bsu{i=1}{r}}
	\newc{\bsim}{\bsu{i=1}{m}}
	
	\newc{\iibr}[2]{\iidx{\bprw{#1}}{#2}}
	\newc{\Intbr}[1]{\iibr{R}{#1}}
	\newc{\intbr}[1]{\iibr{\rg}{#1}}
	\newc{\intt}[3]{\int_{#1}^{#2}\int_\Og~#3~dxdt}
	
	\newc{\itQ}[2]{\dspl{\int\hspace{-2.5mm}\int_{#1}~#2~dz}}
	\newc{\mitQ}[2]{\dspl{\rule[1mm]{4mm}{.3mm}\hspace{-5.3mm}\int\hspace{-2.5mm}\int_{#1}~#2~dz}}
	\newc{\mitQQ}[3]{\dspl{\rule[1mm]{4mm}{.3mm}\hspace{-5.3mm}\int\hspace{-2.5mm}\int_{#1}~#2~#3}}
	
	\newc{\mitx}[2]{\dspl{\rule[1mm]{3mm}{.3mm}\hspace{-4mm}\int_{#1}~#2~dx}}
	\newc{\mitmu}[2]{\dspl{\rule[1mm]{3mm}{.3mm}\hspace{-4mm}\int_{#1}~#2~d\mu}}
	\newc{\iidmu}[2]{\iidx{#1}{#2}}
	
	\newc{\iidm}[3]{{\dspl\int_{#1}~#2~d #3}}
	
	\newc{\itQmu}[2]{\dspl{\int\hspace{-2.5mm}\int_{#1}~#2~d\mu}}
	\newc{\mitQmu}[2]{\dspl{\rule[1mm]{4mm}{.3mm}\hspace{-5.3mm}\int\hspace{-2.5mm}\int_{#1}~#2~d\mu}}

	\newc{\mitQq}[2]{\dspl{\rule[1mm]{4mm}{.3mm}\hspace{-5.3mm}\int\hspace{-2.5mm}\int_{#1}~#2~d\bar{z}}}
	\newc{\itQq}[2]{\dspl{\int\hspace{-2.5mm}\int_{#1}~#2~d\bar{z}}}

	\newc{\pder}[2]{\dspl{\frac{\partial #1}{\partial #2}}}

	\newc{\bdrom}{\bdry{\Og}}

	\newc{\bilhom}{\mbox{Bil}(\mbox{Hom}(\RR^{nm},\RR^{nm}))}
	\newc{\VV}[1]{{V(Q_{#1})}}
	
	\newc{\ccA}{{\mathcal A}}
	\newc{\ccB}{{\mathcal B}}
	\newc{\ccC}{{\mathcal C}}
	\newc{\ccD}{{\mathcal D}}
	\newc{\ccE}{{\mathcal E}}
	\newc{\ccH}{\mathcal{H}}
	\newc{\ccF}{\mathcal{F}}
	\newc{\ccI}{{\mathcal I}}
	\newc{\ccJ}{{\mathcal J}}
	\newc{\ccK}{{\mathcal K}}
	\newc{\ccP}{{\mathcal P}}
	\newc{\ccQ}{{\mathcal Q}}
	\newc{\ccR}{{\mathcal R}}
	\newc{\ccS}{{\mathcal S}}
	\newc{\ccT}{{\mathcal T}}
	\newc{\ccX}{{\mathcal X}}
	\newc{\ccY}{{\mathcal Y}}
	\newc{\ccZ}{{\mathcal Z}}
	
	\newc{\bb}[1]{{\mathbf #1}}
	
	\newc{\myprod}[1]{\langle #1 \rangle}
	\newc{\mypar}[1]{\left( #1 \right)}
	
	\newc{\BLLg}{\mathbf{\LLg}}
	
	\newc{\mA}{\mathbf{A}}
	\newc{\mB}{\mathbf{B}}
	\newc{\mC}{\mathbf{C}}
	\newc{\mD}{\mathbf{D}}
	\newc{\mE}{\mathbf{E}}
	\newc{\mF}{\mathbf{F}}
	\newc{\mJ}{\mathbf{J}}
	\newc{\mG}{\mathbf{G}}
	\newc{\mP}{\mathbf{P}}
	\newc{\mR}{\mathbf{R}}
	\newc{\mQ}{\mathbf{Q}}
	\newc{\mX}{\mathbf{X}}
	\newc{\muu}{\mathbf{u}}
	\newc{\mvv}{\mathbf{v}}
	
	\newc{\mllg}{\mathbb{\lambda}}
	\newc{\mLLg}{\mathbf{\LLg}}

	\newc{\lspn}[2]{\mbox{$\| #1\|_{\Lsp{#2}}$}}
	\newc{\Lpn}[2]{\mbox{$\| #1\|_{#2}$}}
	\newc{\Hn}[1]{\mbox{$\| #1\|_{H^1(\Og)}$}}

	\newc{\mynorm}[2]{\| #1\|_{#2}}

	\newcommand{\RR}{{\rm I\kern -1.6pt{\rm R}}}

	\newc{\itQQ}[2]{\dspl{\int_{#1}#2\,dz}}
	\newc{\mmitQQ}[2]{\dspl{\rule[1mm]{4mm}{.3mm}\hspace{-4.3mm}\int_{#1}~#2~dz}}
	\newc{\MmitQQ}[2]{\dspl{\rule[1mm]{4mm}{.3mm}\hspace{-4.3mm}\int_{#1}~#2~d\mu}}

	\newc{\MUmitQQ}[3]{\dspl{\rule[1mm]{4mm}{.3mm}\hspace{-4.3mm}\int_{#1}~#2~d#3}}
	\newc{\MUitQQ}[3]{\dspl{\int_{#1}~#2~d#3}}

	\newc{\mccP}{\mathbb{P}}
	\newc{\mccK}{\mathbb{K}}
	
	\newc{\DKTmU}{\mccK(U)}
	\newc{\DKTmUold}{(K_U(U)^{-1})^T}
	
	\newc{\myPi}{\mathbf{W}}
	\newc{\myIbar}{\bar{\ccI}_1}
	\newc{\myIhat}{\hat{\ccI}_1}
	\newc{\myIbreve}{\breve{\ccI}_0}
	
	\newc{\mmk}{\mathbf{k}}

	\newcommand{\ma}{\mathbf{a}}
	\newcommand{\mg}{\mathbf{g}}

	\newc{\mfu}{\mathbf{f_u}}
	\newc{\mh}{\mathbf{h}}
	\newc{\mb}{\mathbf{b}}
	\newc{\mf}{\mathbf{f}}

	\newcommand{\barrl}[2]{\barr{ll}\lefteqn{#1}\hspace{#2}&\\}

	\newc{\twomatrix}[1]{\left[\barr{cc}#1\earr\right]}
	\newc{\threematrix}[1]{\left[\barr{ccc}#1\earr\right]}

	\newc{\mN}{\mathbf{N}}
	\newc{\mI}{\mathbf{I}}
	\newc{\mH}{\mathbf{H}}

	\newc{\mk}{\mathbf{k}}
	\newc{\mr}{\mathbf{r}}

	\newc{\DIAGM}[2]{\left[\barr{ccc}#1&0\ldots&0\\
		\vdots&\ddots&\vdots\\0&\ldots0&#2\earr \right]}
	\newc{\DiagM}[2]{\mbox{diag}\left[#1
		\cdots #2 \right]}
	\newc{\vVEC}[2]{\left[\barr{c}#1\\
		\vdots\\#2\earr \right]}
	\newc{\hVEC}[2]{\left[#1
		\cdots #2 \right]}
	
	\newc{\mq}{\mathbf{q}}

	\newc{\msys}[1]{\left\{\barr{l}#1\earr
		\right.}
	\newc{\msysa}[1]{\left\{\barr{ll}#1\earr
		\right.}
	
	\newc{\bbM}{\mathbb{M}}
	\newc{\mat}[1]{\left[\barr{cc}#1\earr\right]}
	
	\newc{\me}{\mathbf{e}}

	\newc{\vecc}[2]{\left[\barr{cc}#1\\#2\earr\right]}
	\newc{\mL}{\mathbb{L}}
	
	\newc{\cO}{{\cal O}}
	
	\newc{\cM}{{\cal M}}
	
	\newc{\myega }{\eg_0(R)}
	\newc{\myeg}{\eg_1(\eg_*)}
	\newc{\myegp}{\hat{\eg}_1(\eg_*)}

\newc{\diagA}{\mathbb{A}_d}
\newc{\mBB}{\mathbb{B}}
\newc{\MLT}[1]{{\cal M}_{lt}(\Og,#1)}
\newc{\ALT}[1]{{\cal A}_{l}(\Og,#1)}
\newc{\mM}{\mathbb{M}}

\newc{\diag}[1]{\mbox{diag}(#1)}
\newc{\off}[1]{\mbox{offdiag}(#1)}
\newc{\mT}{\mathbb{T}}

\usepackage{amssymb}

\newc{\idmu}[2]{{\dspl\int_{#1}~#2~d\mu}}
\newc{\idllg}[2]{{\dspl\int_{#1}~#2~d\llg}}

\begin{document}

	\vspace*{-.8in}
	\begin{center} {\LARGE\em On the smallness of  mean oscillations and regularity of weak solutions to regular/degenerate strongly coupled parabolic systems}
		
	\end{center}

	\vspace{.1in}
	
	\begin{center}

		{\sc Dung Le}{\footnote {Department of Mathematics, University of
				Texas at San
				Antonio, One UTSA Circle, San Antonio, TX 78249. {\tt Email: Dung.Le@utsa.edu}\\
				{\em
					Mathematics Subject Classifications:} 49Q15, 35B65, 42B37.
				\hfil\break\indent {\em Key words:} metric-measure spaces,  H\"older
				regularity, global existence.}}

	\end{center}

	\begin{abstract}
		It will be established that the mean oscillation of bounded weak solutions to strongly coupled parabolic systems is small in small balls. If the systems are regular elliptic then their bounded weak solutions are H\"older continuous. Further assumptions on the systems will even prove that these solutions exist globally. Weak solutions to degenerate systems of porous media type are also studied.  \end{abstract}
	
\section{Introduction} In this paper we investigate bounded weak solutions to the systems
\beqno{exsyszpara0}\left\{\barr{ll}W_t=\Div(\ma(W)DW)+f(W)&\mbox{in $\Og\times (0,T)$,}\\ \mbox{homogeneous BC}&\mbox{on $\partial\Og\times(0,T)$,}\\W(x,0)=W_0(x)&\mbox{on $\Og$,  with $W_0\in W^{1,p_0}(\Og)$ with $p_0>N$} \earr \right.\eeq
on a bounded domain $\Og$ of $\RR^N$, $N\ge2$, with homogeneous Dirichlet or Neumann boundary conditions. Here, $W$ is an unknown vector in $\RR^m$, $m\ge 2$. As usual, $\ma(W), f(W)$ are respectively continuous $m\times m$ matrix and $m\times1$ (Lipschitz) vector. We assume that there are  $\llg(W),\LLg(W)>0$ such that \beqno{elliptic}\llg(W)|\zeta|^2\le \myprod{\ma(W)\zeta,\zeta}\le \LLg(W)|\zeta|^2\quad\forall\zeta\in\RR^m\times\RR^N.\eeq

Under suitable assumptions we will study the following property of weak solutions to \mref{exsyszpara0}

\bdes\item[BMOsmall)] The BMO norm $\|u\|_{BMO(B_R)}$ is small if $R$ is sufficiently small.
\edes
Here, $B_R$ denotes a ball in $\RR^N$ with radius $R$. In \cite{dlebook, dlebook1}, we elucidate the pivotal role played by the smallness of BMO (bounded mean oscillation) norms in small balls within the examination of the regularity of weak solutions (and, in some cases, the global existence of classical solutions) to strongly coupled elliptic/parabolic partial differential systems \cite{Am2, Gius} on $\RR^N$. 

In a recent paper \cite{theo}, we present broader theoretical findings in metric-measure spaces to study BMOsmall) for general functions. In particular, we establish that the mean oscillation of a function on $\RR^{N}$ will be small if $R$ is small, provided there is some basic information on its partial derivatives.  In \refsec{bmosmall-weaksol}, we apply these abstract results  to bounded weak solutions of regular parabolic second-order partial differential systems. In conjunction with the findings in \cite{GiaS, Gius}, we demonstrate that bounded weak solutions exhibit Hölder continuity and hence extend the results in \cite{theo} where we had to assume that  the domain $\Og$ must be sufficiently thin. Assuming further that $\ma$ satisfies a spectral-gap condition, we establish that bounded weak solutions are classical and {\em exist globally}, following the theories of \cite{dlebook, dlebook1}.

We also study the BMOsmall) property of weak solutions to the degenerate system \mref{exsyszpara0} of porous media type in \refsec{degsys}. However, we can not establish the H\"older continuity as it is now well-known for weak solutions to scalar equations (e.g. see \cite{Vasquez}). Finally, we conclude this paper by investigating the uniqueness of weak solutions to the degenerate system \mref{exsyszpara0}  in \refsec{unisec}.

\section{BMOsmall) of bounded weak solutions}\label{bmosmall-weaksol}\eqnoset

In \cite[Theorem 3.1]{theo}, we proved the following  simple theoretic functional result.
\btheo{KEYBMOthm} Let $K\ge1$ and $\Og$ be a domain in $\RR^{N+K}$ and  $\mB_R$ be a cube in $\Og$ with sides parallel to the axes of $\RR^{N+K}$. We write $\RR^{N+K}=\{(x,y_*)\,:\, x\in\RR^N, y_*\in\RR^K\}$.   If for small $R_0,\eg>0$

\bdes
\item [i)] $D_{y_*}W \in L^{p,N+K-p}(\Og)$  for some finite $p\ge1$ and for all $0<R\le R_0$
\beqno{DxNW1}\dspl{\dspl{\int_{B'_R(s_*)}}} \iidx{B_R}{| D_{x_{*}}W(x,y_*)|^p}dy_*\le \eg R^{N+K-p}, \quad B_R\subset\RR^{N},\;B'_R(s_*)\subset\RR^{K};\eeq

\item [ii)]   for {\bf some} $s_*\in\RR^K$  $$f(y_*)=\frac{1}{R^N}\iidx{B_R}{|W(x,y_*)-W_{B_R}|},\quad y_*\in Y.$$
Suppose that $R$ is small such that
$|f(s_*)|\le \eg$ and $|f_{B'_{R}}|\le \eg$.
\edes
Then, for $R,\eg$ as in \mref{DxNW1}, $\mB_R=B'_R\times B_R(s_*)$, and for some constant $C$
\beqno{notBMO}\frac{1}{R^{N+K}}\int_{\mB_R}{|W-W_{\mB_{R}}|}dxdy_*\le C(\eg^\frac{1}{p}+\eg) .\eeq

\etheo

It is important to note that \mref{notBMO} does NOT say that the BMO norm, but the mean oscillation, of a function $W$ on a ball/cube $\mB_R$  with sides parallel to the axes of $\RR^{N+K}$ will be small in $R$ if its mean oscillation in $N$ dimensional ball $B_R$ is also small and \mref{DxNW1} holds for some small $\eg>0$ and $p\ge1$ (a simple information on the $\RR^K$-direction derivatives $D_{x_*}$).

\subsection{Regular strongly coupled parabolic systems} \label{regsys}
We consider the cross diffusion system
\beqno{exsyszpara}\left\{\barr{ll}W_t=\Div(\ma(W)DW)+f(W)&\mbox{in $\Og\times (0,T)$,}\\ \mbox{homogeneous BC}&\mbox{on $\partial\Og\times(0,T)$,}\\W(x,0)=W_0(x)&\mbox{on $\Og$,  with $W_0\in W^{1,p_0}(\Og)$} \earr \right.\eeq
on bounded domain $\Og$ of $\RR^N$ with homogeneous Dirichlet or Neumann boundary conditions. Here, $W$ is an unknown vector in $\RR^m$, $m\ge 2$. As usual, $\ma(W), f(W)$ are respectively continuous $m\times m$ matrix and $m\times1$ (Lipschitz) vector. According to \cite{Am2}, we assume $p_0>N$.

For some positive constants $\llg,\LLg$ we assume the ellipticity condition
\beqno{regellcond}\LLg|DW|^2\ge \myprod{\ma(W)DW,DW}\ge \llg|DW|^2.\eeq

\blemm{BMOreg0} Consider a bounded weak solution $W$ of \mref{exsyszpara}. If the mean oscillation of $W$ on $B'_R$ \beqno{BMOK}\frac{1}{|B'_{R}|}\dspl{\int_{B'_{R}}}|W-W_{B'_R}|^{2}dx_*\le \eg,\eeq then $W$ satisfies i) of \reftheo{KEYBMOthm}. \elemm

\bproof
Let  $B_R$  in \reftheo{KEYBMOthm} be the parabolic cube $Q_R$ in $\RR^N\times(R^2,2R^2)$ then we have $|Q_R|=R^{N+2}$ (for simplicity, we use translations in the $t$-direction). The condition i) of \reftheo{KEYBMOthm} we have to check is
\beqno{DxNW1para}\dspl{\int_{B'_R(s_*)}} \dspl{\int_{R^2}^{2R^2}} \iidx{B_R}{| D_{x_{*}}W(x,y_*)|^p}ds dy_*\le \eg R^{N+2+K-p}, \quad B_R\subset\RR^{N},\;B'_R(s_*)\subset\RR^{K}\eeq
for some $p\ge1$ and $D_{x_{*}}W(x,y_*)$ is the partial derivative of $W$ in the $\RR^K$ direction.
Testing the system with $(W-W_{B'_R})\fg^2\eta$ where $\fg,\eta$ are cut-off function for $B_{2R}, (0,2R^2)$, as usual, we obtain
\beqno{testsys}\dspl{\int_{R^2}^{2R^2}}\dspl{\int_{B'_{R}}}\iidx{ B_{R}}{|DW|^{2}} dx_*ds\le \frac{C}{R^2}\dspl{\int_0^{4R^2}}\dspl{\int_{B_{2R}}}\dspl{\int_{B'_{R}}}|W-W_{B'_R}|^{2}dx_*dxds.\eeq

Therefore, if the mean oscillation of $W$ on $B'_R$ satisfies \mref{BMOK}
$$\frac{1}{|B'_{R}|}\dspl{\int_{B'_{R}}}|W-W_{B'_R}|^{2}dx_*\le \eg,$$ then  we have from \mref{testsys} that
$$\dspl{\int_{R^2}^{2R^2}}\dspl{\int_{B'_{R}}}\iidx{B_{R}}{|DW|^{2}}dx_*ds\le \eg R^{-2} |B'_{R}\times Q_R|.$$

As $|B'_{R}\times Q_R|\sim R^{N+2+K}$, this gives \mref{DxNW1para} when $p=2$. \eproof

As a bounded weak solution $W$ of \mref{exsyszpara} is unique, it  can be approximated by (and so is its mean oscillations by those of ) strong solutions of the linear systems with smooth coefficients
\beqno{approx-exsyszpara}\left\{\barr{ll}W_t=\Div(\ma(W_h)DW)+f(W)&\mbox{in $\Og\times (0,T)$,}\\ \mbox{homogeneous BC}&\mbox{on $\partial\Og\times(0,T)$,}\\W(x,0)=W_0(x)&\mbox{on $\Og$,  with $W_0\in W^{1,p_0}(\Og)$} \earr \right.\eeq
where $W_h$ is the mollifiers of $W$. Therefore, we can talk about the traces of $W$. The smallness of mean oscillations of $W$ on 1-$d$/2-$d$ sections of $\Og$ can be proved by following the argument in \cite{dleJMAA,dlebook1}. In fact, in reducing \mref{approx-exsyszpara} to 2-$d$ sections of $\Og$, we can prove that $\|DW\|_{L^2(\Og\cap\RR^2)}$ is bounded so that, by Poincar\'e's inequality, the smallness of mean oscillations of $W$ follows. Thus, we can verify ii) of \reftheo{KEYBMOthm} and \mref{BMOK}. 

Finally, we can use an induction argument on $N$, $K=1,2$ to assert that

\blemm{BMOreg} Consider a bounded weak solution $W$ of \mref{exsyszpara}. The mean oscillation of $W$ on the parabolic cube $\mQ_R=B'_R\times B_R\times(0,R^2)\subset \RR^{N+K+2}$ is small if $R$ is sufficiently small. \elemm

It is now well known that the higher integrability of the spatial derivatives of $W$ is available under \mref{regellcond} (e.g. see the proof of \cite[Theorem 2.1]{GiaS}, which requires only Caccipoli and Poincar\'e type inequalities in \cite[Lemmas 2.1 and 2.2]{GiaS}). We drop the function $f$ for simplicity and can show that there are $p_*>2$, a constant $C$ such that  $DW\in L^{p_*}_{loc}(\Og\times (0,T))$ such that, via translations, 
\beqno{HIineq}\left(\frac{1}{|\mQ_R|}\dspl{\int_0^{R^2}}\iidx{B_R}{|DW|^{p_*}}ds\right)^\frac{1}{p_*}\le C \left(\frac{1}{|\mQ_{2R}|}\dspl{\int_0^{4R^2}}\iidx{B_{2R}}{|DW|^{2}}ds\right)^\frac{1}{2}, \mbox{ for any $R>0$}.\eeq

This higher integrability of the spatial derivatives of $W$ and BMOsmall) imply the H\"older continuity of $W$ (see \cite[condition (3.2) of Theorem 3.1]{GiaS} using freezing coefficient argument). \reflemm{BMOreg} then yields that
\btheo{Hcontthm} Any bounded weak solution $W$ of \mref{exsyszpara} is H\"older continuous. \etheo

\brem{approxrem} Although the traces of strong solutions are H\"older continuous  on sections of $\Og$ but we can not take their approximation to conclude that the weak solution $u$ is also H\"older continuous as the convergence  (in $L^2$ as in \reflemm{BMOreg}) is too weak. Therefore, we have rely on the theorem of \cite{GiaS} to see that we have the convergence in H\"older spaces in proving \reftheo{Hcontthm}. \erem

On the other hand, suppose that there are $\llg(W),\LLg(W)>0$ such that $$\llg(W)|\zeta|^2\le \myprod{\ma(W)\zeta,\zeta}\le \LLg(W)|\zeta|^2\quad\forall\zeta\in\RR^m\times\RR^N.$$ Let $\nu_*=\sup_{W\in\RR^m}\frac{\llg(W)}{\LLg(W)}$. We say that 
$\ma(W)$ verifies the {\em spectral-gap condition} if  $\nu_*>1-2/N$.  
The spectral gap condition  is just technically sufficient for our theory. In fact, this condition, which is void if $m=1$ (scalar equations) or $N=2$ (planar domains), was used only to guarantee that for some $p>N/2$ and $c>0$ (see \cite{dlebook}) $$\myprod{\ma(W)DX,D(|X|^{2p-2}X)} \ge c\llg(W)|X|^{2p-2}|DX|^2, \quad X\in C^1(\Og,\RR^{m}). $$ 

We then have

\btheo{globalexthm} If $\ma$ satisfies the spectral-gap condition then any bounded weak solution $W$ of \mref{exsyszpara} is classical and exists globally. \etheo

\bproof This is a consequence of \reftheo{Hcontthm} and the theory in \cite{dlebook1}. For the convenience of the readers we sketch here the proof. We follow the proof of \cite[Lemma 4.4.1]{dlebook1} to estimate $\|DW\|_{L^{2p}_{loc}(\Og)}$. Applying the difference operator  $\dg_h$ in $x$ to  the equation of $W$, we see that $W$ weakly solves (dropping the terms involving $\mB$ in the proof of \cite[Lemma 4.4.1]{dlebook1})
$$(\dg_hW)_t=\Div(\ma D(\dg_h W)+\ma_W \myprod{\dg_h W,DW})+f_W  \dg_hW.$$

Let $p\ge1$ and   $\fg,\eta$ be positive $C^1$ cutoff functions for the concentric balls $B_s, B_t$ and the time interval $I$. For any $0<s<t<2R_0$ we test this system with $|\dg_hW|^{2p-2}\dg_hW\fg^2\eta$ and use Young's inequality for the term $|\ma_W||\dg_hW|^{2p-1}|DW||D(\dg_hW)|$ and the spectral gap condition (with $X=\dg_h W$). We get    for some constant $c_0>0$, $\Psi=|f_W|$,  and $Q=\Og_t\times I$
$$\barrl{\sup_{t\in(t_0,T)}\iidx{\Og}{|\dg_hW|^{2p}\fg^2}+c_0\itQ{Q}{ \llg|\dg_hW|^{2p-2}|D(\dg_hW)|^2\fg^2}\le  }{1cm}& C\iidx{\Og\times\{t_0\}}{|\dg_hW|^{2p}\fg^2}+C\itQ{Q}{\frac{|\ma_W|^2}{\llg}|\dg_hW|^{2p}|DW|^2\fg^2}+\\&C\itQ{Q}{[\Psi|\dg_hW|^{2p}\fg^2+|\ma_W|(|\dg_hW|^{2p}|DW|)\fg|D\fg|]}.\earr$$

As $\frac{|\ma_W|^2}{\llg}$ is bounded and $W$ is H\"older continuous (so that $W$ satisfies BMOsmall) on $\Og$), we make use of the weak Gagliardo-Nirenberg inequality with BMO norm \cite[Theorem 2.4.2]{dlebook1}  to estimate the second term on the right hand side, absorb it to the left  and proceed in the same way as in \cite[Lemma 4.4.1]{dlebook1} and obtain the estimate for $\sup_I\|DW\|_{L^{2p}(\Og_R)}$ by letting $h\to0$. As $\sup_I\|DW\|_{L^{2p}(\Og_R)}$ does not blow up for some $p>N/2$, the theory in \cite{Am2} completes the proof. \eproof

\subsection{Degenerate strongly coupled parabolic systems} \label{degsys}
Let   $\mg:\Og\times(0,T)\to\RR^m$ be a bounded map. Consider  the system with Dirichlet or Neumann boundary conditions\beqno{degeqn0} u_t=\Div(\ma(u)Du)+\mg(x,t),\quad u(0)=u_0 \mbox{ on $\Og$, $x\in\RR^{K+2}$}.\eeq

If $\ma(0)=0$, we say that the system is degenerate. For simplicity, we adapt a somewhat strong definition of weak solution $u$ requiring that $\ma(u)Du\in L^2(\Og\times(0,T)$). We have the following result.

\btheo{weakdegthm} Suppose that $\ma(u)$ is positive definite and satisfies the spectral-gap condition, and that $|\ma(u)|$ is bounded. Define the following maps on $\RR^m$
$$\Fg(u)=\int_0^1\ma(su)uds,\; \Psi(u)=\int_0^1\sqrt{\ma(su)}uds.$$ We assume that $\Fg^{-1},\Psi$ are H\"older continuous on $\RR^m$. Then, the bounded weak solution $u$ (which is  assumed to be unique for simplicity) of \mref{degeqn0}  satisfies BMOsmall).

\etheo

It is easy to see that a simple example is $\ma(u) =|u|^\ag A$ for some $\ag>0$ and a constant matrix $A$ satisfying the spectral-gap condition.
The BMOsmall) property for solutions to parabolic systems in \reftheo{weakdegthm} and \reflemm{BMOreg0} are referring to the parabolic cube $B'_R\times B_R\times (t_0,t_0+R^2)$. However, these assertions can be improved (see \refrem{paracuberem} after the proof).

\bproof
Consider the strong solutions $u_{\eg}$ of the following nonlinear systems 
\beqno{approx-deg}\left\{\barr{ll}u_t=\Div(\ma^{(\eg)}(u)Du)+g(x,t)&\mbox{in $\Og\times (0,T)$,}\\ \mbox{homogeneous BC}&\mbox{on $\partial\Og\times(0,T)$,}\\u(x,0)=u_0(x)&\mbox{on $\Og$,  with $u_0\in W^{1,p_0}(\Og)$} \earr \right.\eeq
Here, $\ma^{(\eg)}(u)=\ma(u)+\eg Id$ for $\eg>0$. Because $|\ma(u)|$ is bounded and $\ma(u)$ satisfies the spectral-gap condition, the above system is regular elliptic and $u_\eg$ exists on $\Og\times(0,T)$ (see \reftheo{globalexthm}), we can talk about the traces of the strong solutions $u_\eg$ below.

First, fixing $x_*$ (reducing the system to 2-$d$ sections of $\Og$) we consider the following planar system $$u_t=\Div_x(\ma^{(\eg)}(u)D_xu)+\mg(x,x_*,t), \quad x\in\Og\cap(\{x_*\}\times\RR^2), t\in(0,T).$$

We test this system with $\ma^{(\eg)}(u_\eg) (u_{\eg})_t$ as in \cite{dleJMAA}, $u_{\eg}$ is classical, to see that $\|\ma^{(\eg)}(u_\eg)D_x(u_{\eg})\|_{L^2(\Og\times\RR^2)}$, the norm on 2-$d$ sections of $\Og$,  is uniformly bounded in $\eg$ (note that the proof in \cite{dleJMAA}  requires only that $\ma^{(\eg)}$ is regular elliptic and $|\ma^{(\eg)}(u_\eg)|$ is bounded as we can see in \cite[Section 4.2.2.1]{dlebook1}). Because $\Fg_u(u)=\ma(u)$ is positive definite, from the definition of $\ma^{(\eg)}$, we see that 
$$|D_x\Fg(u_{\eg})|^2\le \eg^2|D_x(u_\eg)|^2 +2\eg\myprod{\Fg_u(u_\eg)D_xu_{\eg},D_xu_{\eg}}+|\Fg_u(u_\eg)D_xu_{\eg}|^2 =|\ma^{(\eg)}(u_\eg)D_x(u_{\eg})|^2$$
so that 
$\|D_x\Fg(u_{\eg})\|_{L^2(\Og\times\RR^2)}$ is uniformly bounded in $\eg$ also. Thus,  $\Fg(u_{\eg})$ satisfies BMOsmall) uniformly by Poincar\'e's inequality. Let $\eg\to0$, the corresponding sequence   $\{\Fg(u_{\eg})\}$ converges to  $\Fg(u)$ in $L^2(\Og)$ (with $u$ is the weak solution of \mref{degeqn0}, see also \refrem{uniquerem}). So that $\Fg(u)$ satisfies BMOsmall).  As $\Fg^{-1}$ is H\"older continuous, $u$ also verifies BMOsmall) on $\Og\times (0,T)$.  

Denote by $B'_R, B_R$ be balls of radius $R$ in $\RR^2, \RR^K$ respectively. We start our induction on $K$ with $K=1,2$.

Setting $W=\Psi(u)=\int_0^1\sqrt{\ma(su)}uds$ (note that $\Fg_u(u)=\ma(u)$), we have $DW=\sqrt{\Fg_u(u)}Du$. Testing the system with $(u-u_{B'_R})\fg^2\eta$ where $\fg,\eta$ are cut-off function for $B_{2R}, (R^2,2R^2)$, we obtain as usual (compare to \mref{testsys} as $\myprod{D\Fg(u),Du}=|DW|^2$) 
\beqno{testsys1}\dspl{\int_{R^2}^{2R^2}}\dspl{\int_{B'_{R}}}\iidx{ B_{R}}{|DW|^{2}} dx_*ds\le \frac{C}{R^2}\dspl{\int_{R^2}^{2R^2}}\dspl{\int_{B_{2R}}}\dspl{\int_{B'_{R}}}|u-u_{B'_R}|^{2}dx_*dxds.\eeq
As $u$ also satisfies BMOsmall) for the ball $B'_R\subset\RR^K$, $K=1,2$, we derive that
$$\dspl{\int_{R^2}^{2R^2}}\dspl{\int_{B'_{R}}}\iidx{B_{R}}{|DW|^{2}}dx_*ds\le \eg R^{-2} |B'_{R}\times Q_R|,  $$
where $ Q_R$ is the parabolic cube $B_R\times(R^2,2R^2)$. As $|B'_{R}\times Q_R|\sim R^{N+2+K}$, this gives i) of \reftheo{KEYBMOthm} when $p=2$ for $N=2, K=1,2$.

If $\Psi$ is H\"older continuous we can apply \reftheo{KEYBMOthm} to $W=\Psi(u)$ and its ii)  on $\Og\cap\RR^N\times (0,T)$) is verified (for $N=2$ and because $u$ satisfies BMOsmall)). We then conclude that $u=\Psi^{-1}(W)$ satisfies BMOsmall) on $\Og\times (0,T)$.

Now the proof follows easily by induction with $N=3,4,\ldots$ and $K=1,2$. \eproof

\brem{uniquerem} The assertion also holds if $g$ in \mref{degeqn0} is replaced by $g(x,t)u$ for some matrix $g(x,t)$. If there is no uniqueness for weak solutions then we can assert that there is a weak solution which satisfies BMOsmall). The proof is similar, we keep working with strong solutions of \mref{approx-deg} and take to the limit $\eg\to0$ at the end to obtain such a weak solution. Note that in the proof, we can talk about $D\Fg(u), D\Psi(u)$ as functions in $L^2(\Og)$, $\Og\subset\RR^2$, but not $Du$ (which exists in the distribution sense but not a function in $L^2(\Og)$). Concerning the uniqueness, we refer the readers to \refsec{unisec}.

\erem

\brem{book1rem} We should also remark here that although the assertion of \cite[Corollary 4.2.7]{dlebook1} is valid but its proof (using scalar equations techniques) can be replaced by that of \cite[Lemma 3.7.1]{dlebook1} where we already proved that $\|Du\|_{L^{2p}(B_R)}$ is bounded for some $p>1$ provided that $u$ satisfies BMOsmall) on $\Og$ (the spectral-gap condition is void  when $\ma$ is regular and $N=2$, see also the proof of \reftheo{globalexthm}). \erem

\brem{paracuberem} We should emphasize that the BMOsmall) property for solutions to parabolic systems in \reftheo{weakdegthm} and \reflemm{BMOreg0} are referring to the {\em parabolic cube} $B'_R\times B_R\times (t_0,t_0+R^2).$  In fact, like \refrem{approxrem}, we have the convergence of strong solutions of \mref{approx-deg} (or \mref{approx-exsyszpara}) in $L^2(\Og)$ to see that $u$ satisfies 
BMOsmall) on $\Og$. Indeed, $|D\Fg(u_\eg)|^2$'s are uniformly integrable  on 2-$d$ dimensional sections of $\Og$ and we can invoke Fubini's theorem to conclude that $|D\Fg(u_\eg)|^2$'s are uniformly integrable on $\Og$ as well and therefore $\{\Fg(u_\eg)\}$ is compact in $L^2(\Og)$. The strong solutions are H\"older continuous but we can not take to the limit in $L^2(\Og)$ to say the same about $u$. However, as $|D\Fg(u_\eg)|^2$'s are uniformly integrable  on  $\Og$, $\Fg(u_\eg)$'s satisfy BMOsmall) on $\Og$ uniformly and so do $u_\eg$'s.  The mean oscillations of $u_\eg$'s converge to that of $u$. Thus, the BMOsmall) property for solutions to parabolic systems in \reftheo{weakdegthm} and \reflemm{BMOreg0} are referring to the {\em elliptic ball}  $B'_R\times B_R$ also.

\erem

Although, in the degenerate case \mref{degeqn0}, a higher integrability for $D\Fg(u)$ in $\Og\times(0,T)$ can be established by mimicking the proof in \cite{GiaS}. But using this together with BMOsmall), we could not prove the H\"older continuity of $u$ (see \reftheo{Hcontthm}) by comparing it with solutions of constant coefficient systems as in the regular case, even we assume $\ma(u)=\Fg_u(u)$ satisfies the spectral-gap condition. Howerver, if \mref{degeqn0} is a scalar equation then the H\"older continuity of $u$ is well known, see \cite{Vasquez}.

	\section{Uniqueness of weak solutions} \label{unisec}\eqnoset  First of all we adapt the following concept of weak solutions to the equation
$$\frac{d}{dt}u=\Div(D\Fg(u))+F(x,t)\mbox{ in $\Og\times(0,T)$},\quad u(0)=u_0 \mbox{ on $\Og$}.$$

We say that $u$ is a weak solution on $Q_T=\Og\times(0,T)$ if 
\bdes
\item [i)] $u\in L^1(Q_T)$, $\Fg(u)\in L^1(0,T:W^{1,1}_0(\Og))$;
\item [ii)] for all test function  $\eta\in C^1(\bar{Q_T})$, $\eta=0$ on the parabolic boundary of $Q_T$ we have ($dz=dxdt$)
$$\itQ{Q_T}{[D\Fg(u)D\eta-u\eta_t]}=\itQ{Q_T}{F\eta}+\iidx{\Og}{u_0\eta(x,0)}.$$

\edes

If $F$ is independent of $u$ then it is well known that there is a weak solution and it is unique (see \cite[section 5.3]{Vasquez}). We can generalize the uniqueness part of  this result if we take the following  stronger definition
\bdes \item [W)] $u(\cdot,t), D\Fg(u(\cdot,t))\in L^2(\Og)$ for all $t\in(0,T)$ and we can use $u$ as a test function in ii).
\edes

One can combine Steklov's average and mollifiers to see that W) is reasonable. We then have the following result.

\btheo{uniquethm} Let $\Fg(u)$ be a function on $\RR$, $\mg$ is a function on $\Og\times(0,T)$. Consider  the equation \beqno{degeqn} \frac{d}{dt}u=\Div(D\Fg(u))+\mg(x,t)u,\quad u(0)=u_0 \mbox{ on $\Og$}.\eeq
Then a weak solution $u$ satisfying W) of this equation is  unique  if \beqno{unicond}\Fg_u(0)=0,\; \Fg_{uu}(tv)v\ge0 \mbox{ for $t\in(0,1)$ and $v\in\RR$},\; \dspl{\int_0^T}\iidx{\Og}{\mg(x,s)}ds<\infty.\eeq

\etheo

The (local) existence of a weak solution of \mref{degeqn} can be established by several means, among these is the approximation method.
A typical example of \mref{degeqn} with condition \mref{unicond} is: Let $\Fg(u)=|u|^ku$ be a function, $k\ge1$. We have $D\Fg(u)=(k+1)|u|^{k-1}Du$ and $\Fg_{uu}(tv)v=(k+1)k|t|^{k-2}t|v|^{k}\ge0$ for $t\in(0,1)$. Then the theorem is applied. Also, we can remove W) if one make use of Steklov's average and mollifiers (but we still have to assume $Du(\cdot,t)\in L^2(\Og)$ for all $t\in(0,T)$). For simplicity, we assume W) here.

Before presenting the proof of \reftheo{uniquethm} let us recall the following well known Gr\"onwall inequality
\blemm{gronwallineq} Let $y(t), q(t)$ be functions on $(0,T)$. Assume that $y$ is a.e. differentiable and $q(t)$ is integrable. Assume that $$\frac{ d}{dt}y(t)\le q(t)y(t)+c,\; y(0)=y_0\quad\mbox{a.e. on $(0,T)$}, $$
where $c$ is a constant. Then $y(t)\le e^{\int_0^t q(s)ds}[y_0+c]$.

\elemm

{\bf Proofof \reftheo{uniquethm}:} For any functions $a,b$ (assuming sufficient integrabilities) we have 
$$\myprod{D(\Fg(a)-\Fg(b)),D(a-b)}=\int_0^1\myprod{\Fg_u(sa+(1-s)b)D(a-b), D(a-b)} ds,$$
$$\Fg_u(sa+(1-s)b)=\int_0^1 \Fg_{uu}(t(sa+(1-s)b))(sa+(1-s)b)dt+\Fg_u(0).$$
So that, as $\Fg_u(0)=0$,
$$\myprod{D(\Fg(a)-\Fg(b)),D(a-b)}=\int_0^1\int_0^1\myprod{ \Fg_{uu}(t(sa+(1-s)b))(sa+(1-s)b)D(a-b), D(a-b)}dt ds.
$$
We see that  $\myprod{D(\Fg(a)-\Fg(b)),D(a-b)}\ge0$,as $\Fg_{uu}(tv)v\ge0$ for $t\in(0,1)$ with $v=(sa+(1-s)b)$.

Consider two weak solutions $u,v$ to the equation \mref{degeqn}. Subtracting, we see that
$w=u-v$ satisfies
$\frac{d}{dt}w=\Div(D(\Fg(u)-\Fg(v)))+\mg(x,t)w$. If $u,v$ satisfy W), we can test this equation with  $w$ (using Steklov's average and mollifiers) to obtain
$$\frac{d}{dt}\iidx{\Og}{|w|^2}+\iidx{\Og}{\myprod{D(\Fg(u)-\Fg(v)),D(u-v)}}= \iidx{\Og}{\mg(x,t)|w|^2}.$$
As $\myprod{D(\Fg(u)-\Fg(v)),D(u-v)}\ge0$, we can drop the nonegative term on the left hand side to have
$$\frac{d}{dt}\iidx{\Og}{|w|^2}\le \iidx{\Og}{\mg(x,t)|w|^2},\; w(0)=0$$
so that by Gr\"onwall 's inequality \reflemm{gronwallineq} and $y_0=c=0$, as  $\dspl{\int_0^T}\iidx{\Og}{\mg(x,s)}ds<\infty$, we have $w\equiv0$ on $\Og\times(0,T)$. \eproof

\brem{unisysrem} The above calculation can be extended to the system like \mref{degeqn} when $\Fg:\RR^m\to\RR^m$ with $\Fg_{uu}(tv)v$ is a nonnegative definite matrix for $t\in(0,1)$, $v\in\RR^m$. 
Obviously, we can have $\mg(x,t)=\mbox{diag}[\mg_1(x,t),\ldots,\mg_m(x,t)]$ with $\dspl{\int_0^T}\iidx{\Og}{\mg_i(x,s)}ds<\infty$.

\erem

\brem{degeqnrem} The systems like \mref{degeqn} are also written as $\frac{d}{dt}u=\Div(\hat{\Fg}(u)Du)+q(x,t)u$, with $\hat{\Fg}(u)$ being a $m\times m$ matrix.
Then a weak solution of this system is unique  if $\hat{\Fg}(0)=0$, $\hat{\Fg}_{u}(tv)v$ is a nonnegative definite matrix for $t\in(0,1)$ and $v\in\RR^m$. The same conditions apply for $\mg_i$'s and we can use approximation to obtain H\"older continuity for weak solutions satisfying W).
\erem

\bibliographystyle{plain}

\begin{thebibliography}{10}
	
	
	\bibitem{Am2} H. Amann, \newblock{Dynamic theory of quasilinear parabolic systems III. Global existence,} {\em  Math Z.} 202 (1989), pp. 219–-250.
	
	
	\bibitem{GiaS} M. Giaquinta and M. Struwe, \newblock{ On the partial regularity of weak solutions of nonlinear parabolic
		systems}. {\em  Math. Z.}, 179(1982),  437--451.
	
	
	\bibitem{Gius}
	E. Giusti.
	\newblock {\em Direct Methods in the Calculus of Variations}.
	\newblock World Scientific, Singapore, 2003.
	
	
	\bibitem{dlebook} D. Le, \newblock{\em Strongly Coupled Parabolic and Elliptic Systems: Existence and Regularity of Strong/Weak Solutions.} De Gruyter, 2018.
	
	\bibitem{dlebook1} D. Le, \newblock{\em Cross Diffusion Systems: Dynamics, Coexistence and Persistence.} De Gruyter, 2022.
	
	
	\bibitem{dleJMAA} D. Le, \newblock{On the global existence of a generalized Shigesada-Kawasaki-Teramoto system,} {\em J. Math. Anal. App.} 2021.
	
	\bibitem{theo} D. Le, \newblock{On the smallness of mean oscillations on metric-measure spaces and applications
	}, {\em preprint} arXiv:2403.06696.
	
	\bibitem{Vasquez} J. L. Vasquez, {\em The Porous Medium Equation
		Mathematical Theory}, Oxford-CLARENDON PRESS 2007.
	
	
\end{thebibliography}

\end{document}